\documentclass[12pt,a4paper,leqno]{amsart}
\usepackage{amsmath}
\usepackage{amssymb}

\newtheorem{theorem}{Theorem}[section]

\newtheorem*{theoremA}{Theorem  A} 
\newtheorem*{theoremB}{Theorem  B} 
\newtheorem*{corollaryC}{Corollary C} 
 
\newtheorem{lemma}[theorem]{Lemma}
\newtheorem{proposition}[theorem]{Proposition}

\newcommand{\D}{{\bf D}}        
\newcommand{\R}{{\bf R}}       
\newcommand{\Z}{{\bf Z}}       

\catcode`\@=11 \@addtoreset{equation}{section} \catcode`\@=12

\title[ Trace inequalities] {On $L^p$--$L^q$ trace inequalities}
\author[Cascante]{Carme Cascante}
\address{Departament de Matem\`atica Aplicada i An\`alisi,
Facultat de Matem\`atiques,  Universitat de Barcelona, Gran Via
585, 08071~Barcelona, Spain} \email{cascante@ub.edu}

\author[Ortega]{Joaquin M. Ortega}
\address{ Departament de Matem\`atica Aplicada i An\`alisi,
Facultat de Matem\`atiques,   Universitat de Barcelona, Gran Via
585, 08071~Barcelona, Spain} \email{ortega@ub.edu}

\author[Verbitsky]{Igor E. Verbitsky}
\address{ Department of Mathematics,  University of Missouri-Columbia, 
Columbia,  MO 65211, USA} \email{igor@math.missouri.edu}
\thanks{First two authors partially supported by DGICYT Grant
PB98-1242-C02-01, and CIRIT Grant 2001SGR00172.}

\date{}
\begin{document}
\begin{abstract} 
We  give necessary and sufficient conditions in order that inequalities of the type  
$$
\| T_K f\|_{L^q(d\mu)}\leq C \, \|f\|_{L^p(d\sigma)}, \qquad f \in
L^p(d\sigma),
$$  
hold for a class of integral operators $T_K f(x) = \int_{R^n} K(x, y) \, f(y) \, d \sigma(y)$ with nonnegative 
kernels,  
and measures $d \mu$ and $d\sigma$ on $\R^n$, in the case where 
$p>q>0$ and $p>1$.    

An important model  is provided by the 
dyadic integral operator with  kernel $K_{\mathcal D}(x, y) \sum_{Q\in{\mathcal D}}  K(Q) \chi_Q(x) \, \chi_Q(y)$, where $\mathcal
D=\{Q\}$ is the family of  all dyadic  cubes in $\R^n$, 
and $K(Q)$ are arbitrary nonnegative constants associated with $Q \in{\mathcal D}$. 

The corresponding continuous versions are 
deduced from their dyadic counterparts. In particular, we show that, for the convolution operator 
$T_k f = k\star f$ with positive radially decreasing kernel $k(|x-y|)$, the trace inequality 
$$
\| T_k f\|_{L^q(d\mu)}\leq C \, \|f\|_{L^p(d x)}, \qquad f \in
L^p(dx),
$$ 
holds  if and only if  ${\mathcal W}_{k}[\mu] \in L^s (d\mu)$, where 
$s = {\frac{q(p-1)}{p-q}}$. Here  ${\mathcal W}_{k}[\mu]$ is a nonlinear Wolff 
 potential  defined by 
 $
 {\mathcal W}_{k}[\mu](x)=\int_0^{+\infty}  k(r)
\overline{k}(r)^{\frac 1 {p-1}} \mu (B(x,r))^{\frac 1{p-1}} \, r^{n-1} \, dr, 
 $
 and 
 $
\overline{k}(r)=\frac1{r^n}\int_0^r
k(t) \, t^{n-1} \, dt$. Analogous inequalities for $1\le q < p$ were 
characterized  earlier by the authors using a different method which is not applicable 
when $q<1$.

 \end{abstract} \maketitle

\section{Introduction}  \label{introduction}

 We  consider  inequalities of the
type
 $$
\| T_K f\|_{L^q(d\mu)}\leq C \, \|f\|_{L^p(d\sigma)}, \qquad f \in
L^p(d\sigma),
$$  
  where 
$d\sigma$ and $d\mu$  are  nonnegative  Borel measures  on $\R^n$,
and
$$T_K [f  d \sigma] (x) = \int_{\R^n} K(x, y) \, f(y) \, d \sigma (y), \qquad x \in
\R^n,$$ is an integral operator with nonnegative kernel $K(x, y)$. 
Our goal is to obtain explicit necessary 
and sufficient conditions on $K$, $\mu$, and $\sigma$ which ensure that such inequalities hold 
for a broad class of  dyadic and radially nonincreasing kernels in the case $0<q<p<\infty$ and $p>1$.  
 
 We observe that  such trace inequalities for $q=1$ are related to 
the so-called Wolff inequality and nonlinear Wolff potentials which appeared originally in 
\cite{hedbergwolff} in connection 
with deep approximation properties of  Sobolev spaces. (See \cite{adamshedberg} 
where nonlinear potential theory  and its applications are presented.) Generalizations 
of Wolff's inequality and appropriate versions of Wolff's potentials  served as our main tools in 
earlier work 
\cite{cascanteortegaverbitsky}, \cite{cascanteortegaverbitsky2} where trace inequalities were characterized 
in the ``upper triangle case'' $1\le q<p<\infty$. However, methods used   
there where based on duality arguments which are not applicable if $0<q<1$. 
 
 In this paper, we are mostly interested in the difficult case where $0<q<1$ and $p>1$. 
  Notice that our approach also works for 
 $1\le q<p<\infty$. 
We employ dyadic models, two weight Carleson measure theorems, and certain imbedding and 
 interpolation theorems for discrete Littlewood--Paley spaces. 
Some of these techniques developed in the framework of trace inequalities turned out to be useful in applications to quasilinear and fully nonlinear PDE  
 (see \cite{phucverbitsky}).

An important special case is given by  dyadic  integral operators
$T_{K_{\mathcal D}}$ and the corresponding integral inequalities from which continuous versions 
will be deduced. Let $ {\mathcal D}=\{ Q\}$ be the family of all dyadic cells
 $Q = 2^{i} (k + [0, \, 1)^n)$,  $i \in \Z, \, k \in \Z^n$, and let $r_Q$ denote 
 the side length of $Q$. For $K: \, {\mathcal D}\rightarrow \R^+$, the kernel $K_{\mathcal D}(x,y)$ on
$\R^n\times\R^n$ is defined by  
$$K_{\mathcal D}(x,y) = \sum_{Q\in{\mathcal
D}}  K(Q) \,  \chi_Q(x) \, \chi_Q(y),$$ 
where $\chi_Q$ is the characteristic function of $Q \in {\mathcal
D}$.

Let $\nu$ be a nonnegative Borel measure on $\R^n$, and let $f \in L^1_{\text{loc}} (d \nu)$. 
We define the  dyadic integral 
operator:
  $$ T_{K_{\mathcal D}}[f \, d\nu] (x)=  \int_{\R^n} K_{\mathcal D} (x, y) f(y) \,
d\nu(y)=\sum_{Q\in{\mathcal D}}  K(Q)
\chi_Q(x) \int_Q f \, d \nu.$$   
We denote by $\overline{K}(Q)(x)$
the function  
$$ \overline{K}(Q)(x)=\frac1{\sigma(Q)}\sum_{Q'\subset Q} K(Q')
\sigma(Q')\chi_{Q'}(x),  
$$
supported on $Q$, where $\sigma$ is a nonnegative Borel measure on $\R^n$, and 
$ \overline{K}(Q)(x)=0$ if $\sigma(Q)=0$. 

 For $x\in\R^n$, the generalized Wolff potential was defined in \cite{cascanteortegaverbitsky2} 
as  
\begin{equation}\label{wolff1}
{\mathcal{W}}_{K, \,\sigma}^{\mathcal D} [\nu](x)= \sum_{Q\in{\mathcal D}} K(Q) 
\left[ \int_{ Q} \overline{K}(Q)(y) \, d\nu(y) \right]^{p'-1}  \sigma(Q) \, \chi_Q(x). 
\end{equation}

Let us assume that the pair $(K, \sigma)$  satisfies  the dyadic logarithmic bounded
oscillation condition {\rm (\text{DLBO})}:
  \begin{equation}\label{dlbocondition}
  \sup_{x \in Q} \overline{K}(Q)(x) \le A \, \inf_{x \in Q}
\overline{K}(Q)(x),
 \end{equation}
 where $A$ does not depend on $Q \in {\mathcal D}$.  If $K$ is a radially nonincreasing kernel  and $d\sigma=dx$, or if $K(Q)=r_Q^{n-\alpha}$, $0<\alpha<n$, and $\sigma$ satisfies an adequate dyadic reverse doubling condition, then the pair $(K,\sigma)$ satisfies condition {\rm (\text{DLBO})} (see \cite{cascanteortegaverbitsky2}).

For $(K, \sigma) \in {\rm (\text{DLBO})}$, we set $\overline{K}(Q) = \inf_{x \in Q}\overline{K}(Q)(x)$, 
$Q \in {\mathcal D}.$ In this case the definition of Wolff's potential can be simplified:    
\begin{equation}\label{wolffpotential}
{\mathcal{W}}_{K, \,\sigma}^{\mathcal D} [\nu](x)=\sum_{Q\in{\mathcal
D}} K(Q) \,  [\overline{K}(Q)]^{p'-1} \, [\nu(Q)]^{p'-1}  \sigma(Q) \, \chi_Q(x). 
\end{equation}

 The following theorem is proved in Section~\ref{section2}. 
  
\begin{theoremA}\label{theoremAA}
Let $K:{\mathcal D}\rightarrow\R^+$, $0<q<p<+\infty$, $1<p<+\infty$.  Let $\mu$ and $\sigma$ be 
 nonnegative Borel measures on $\R^n$. Suppose that  $(K, \sigma) \in   {\rm (DLBO)}$. 
 Then the trace inequality
\begin{equation}\label{00.10}
\int_{\R^n} \left \vert \, T_{K_{\mathcal D}}[fd\sigma] \, \right \vert^q 
\,d\mu \leq  C \, \|f\|_{L^p(d\sigma)}^q,  \qquad f\in
L^p(d\sigma),
 \end{equation}
 holds if and only if 
${\mathcal W}_{K, \,\sigma}^{\mathcal D} [\mu]\in
L^{\frac{q(p-1)}{p-q}}(d\mu).$ 
\end{theoremA}

In Section \ref{section3}, we treat continuous trace inequalities.
We state a  version of Theorem~A for
integral operators with radial  kernels, 
$$T_k[f \, d \sigma](x) = \int_{\R^n} k(|x-y|) \, f(y) \, d \sigma(y).$$
Here 
$k=k(r)$, $r>0$, is an arbitrary lower
semicontinuous  nonincreasing positive function.  The corresponding nonlinear potential is defined by  
$$
 {\mathcal W}_{k,\,\sigma}[\mu](x)=\int_0^{+\infty} k(r) \, 
\sigma(B(x,r)) \, 
\left(\int_{B(x,r)}\overline{k}(r)(y) \, d\mu(y)\right)^{p'-1}  \,
\frac{dr}{r},  $$
where 
$$
\overline{k}(r)(x)=\frac1{\sigma(B(x,r))}\int_0^r
k(s) \, \sigma(B(x,s)) \, \frac{ds}{s},$$ 
for $x\in\R^n$, $r>0$.

\begin{theoremB}\label{theoremBB}
 Let $0< q<p<+\infty$, $1<p< \infty$, 
 and let $\mu$ and $\sigma$ be nonnegative Borel
measures on $\R^n$. Assume that  $\sigma$ satisfies a  doubling condition, and   the pair $(k,
\sigma)$  has  the  following logarithmic bounded oscillation property
{\rm (\text{LBO})}:  
 \begin{equation}\label{lbocondition}
 \sup_{y \in B(x,r)} \overline{k}(r)(y) \le A \, \inf_{y \in B(x,r)} \overline{k}(r)(y),
 \end{equation}
 where $A$ does not depend on
$x\in\R^n, \, r>0$.  Then the trace inequality
\begin{equation}\label{tracecont}
\int_{\R^n} \left \vert \, T_{k}[fd\sigma] \, \right \vert^q 
\,d\mu \leq  C \, \|f\|_{L^p(d\sigma)}^q,  \qquad f\in
L^p(d\sigma),
\end{equation}
 holds if and only if 
${\mathcal W}_{k, \,\sigma} [\mu]\in
L^{\frac{q(p-1)}{p-q}}(d\mu)$.
  \end{theoremB}

 The  {\rm (\text{LBO})} property  is satisfied
  by all radially nondecreasing kernels in the case $d \sigma = dx$, or by Riesz kernels 
$k(x) = |x|^{\alpha-n}$, $0<\alpha<n$, if $\sigma$ satisfies an adequate reverse doubling condition (see \cite{cascanteortegaverbitsky2}). In particular, the following corollary holds for convolution 
operators $T_k[f] = k\star f$ and $d \sigma =dx$. 

 \begin{corollaryC}\label{convolution} Let $0< q<p<+\infty$ and $1<p< \infty$. Let $\mu$  be a 
 nonnegative Borel
measure on $\R^n$. Suppose $k=k(|x-y|)$, where $k(r)$ is  a 
lower semicontinuous nonincreasing positive function on $\R_+$, 
and  $
\overline{k}(r)=\frac1{r^n}\int_0^r
k(t) \, t^{n-1} \, dt$. 
Then the trace inequality 
\begin{equation}\label{traceconv}
||k\star f||_{L^q(d \mu)} \leq  C \, \|f\|_{L^p(dx)},  \qquad f\in
L^p(dx),
\end{equation} 
 holds if and only if 
 \begin{equation}\label{wolffpotentialc}
 {\mathcal W}_{k}[\mu](x)=\int_0^{+\infty}  k(r) \, 
\overline{k}(r)^{\frac 1{p-1}} \mu (B(x,r))^{\frac 1{p-1}} \, r^{n-1} \, dr \in L^{{\frac{q(p-1)}{p-q}}}(\mu).
\end{equation}
 \end{corollaryC} 
 \bigskip
   
 \noindent {\bf Remark 1.}  For Riesz kernels, 
  a proof of Corollary~C was sketched previously in \cite{verbitsky2}. However, some 
  technical details 
  related to passing from discrete to continuous versions using shifts of the dyadic 
  lattice are given below for the first time.  \medskip
  
   \noindent {\bf Remark 2.} A more complicated characterization of (\ref{traceconv})  for Bessel kernels 
    was given earlier in \cite{mazyanetrusov} in terms of a certain capacity condition.

\section{Discrete trace inequalities}\label{section2}

In this section we will consider discrete trace inequalities. Recall that if
$ {\mathcal D}=\{ Q\}$ denote
 the family of all dyadic cells, 
and $K: \, {\mathcal D}\rightarrow \R^+$, we have defined the kernel $K_{\mathcal D}(x,y)$ on
$\R^n\times\R^n$  by  $$K_{\mathcal D}(x,y) = \sum_{Q\in{\mathcal
D}}  K(Q) \,  \chi_Q(x) \, \chi_Q(y).$$

If $\nu$ is a nonnegative Borel measure on $\R^n$, and  $f \in L^1_{\text{loc}} (d \nu)$, 
 the  dyadic integral 
operator is given by 
  $$ T_{K_{\mathcal D}}[f \, d\nu] (x)=  \int_{\R^n} K_{\mathcal D} (x, y) f(y) \,
d\nu(y)=\sum_{Q\in{\mathcal D}}  K(Q)
\chi_Q(x) \int_Q f \, d \nu.$$   
 If $f\equiv 1$, we write $ T_{K_{\mathcal D}}[\nu]$ in place of $ T_{K_{\mathcal D}}[f \, d\nu]$.

Let $\sigma$ be another nonnegative Borel measure on $\R^n$. 
As was already mentioned in the introduction, we will assume that the pair $(K, \sigma)$  
satisfies  the dyadic logarithmic bounded
oscillation condition {\rm (\text{DLBO})}:
  \begin{equation}\label{dlbo}
  \sup_{x \in Q} \overline{K}(Q)(x) \le A \, \inf_{x \in Q}
\overline{K}(Q)(x),
 \end{equation}
 where $A$ does not depend on $Q \in {\mathcal D}$, and $$ \overline{K}(Q)(x)=\frac1{\sigma(Q)}\sum_{Q'\subset Q} K(Q')
\sigma(Q')\chi_{Q'}(x).  
$$

For $(K, \sigma) \in {\rm (\text{DLBO})}$, we set $\overline{K}(Q) = \inf_{x \in Q}\overline{K}(Q)(x)$, 
$Q \in {\mathcal D}$, where  $\overline{K}(Q)=0$ if $\sigma(Q)=0$. 
Then    the generalized Wolff potential  introduced in \cite{cascanteortegaverbitsky2} 
can be defined alternatively in an equivalent way by:   
\begin{equation}\label{wolff2}
{\mathcal{W}}_{K, \,\sigma}^{\mathcal D} [\nu](x)=\sum_{Q\in{\mathcal
D}} K(Q) \,  [\overline{K}(Q)]^{p'-1} \, [\nu(Q)]^{p'-1}  \sigma(Q) \, \chi_Q(x). 
\end{equation}
The generalized Wolff inequality proved in \cite{cascanteortegaverbitsky2} (Theorem A) yields: 
\begin{equation}\label{wolffinequality}
C_1 \,  \int_{\R^n} {\mathcal{W}}_{K, \,\sigma}^{\mathcal D} [\nu] \, d \nu \le 
 \int_{\R^n}  \left (T_{K_{\mathcal D}}[\nu]\right )^{p'} \,d \sigma 
  \le C_2 \,   \int_{\R^n} {\mathcal{W}}_{K, \,\sigma}^{\mathcal D} [\nu] \, d \nu,
 \end{equation}
 i.e., the energy of the measure $\nu$ given by $\int_{\R^n} 
 \left (T_{K_{\mathcal D}}[\nu]\right )^{p'} d \sigma$ 
 is equivalent, 
under the {\rm (\text{DLBO})} assumption,    to 
  $$ \int_{\R^n} {\mathcal{W}}_{K, \,\sigma}^{\mathcal D} [\nu] \, d \nu = 
\sum_{Q\in{\mathcal
D}} K(Q) \, \sigma(Q) [\overline{K}(Q)]^{p'-1} [\nu(Q)]^{p'}.
$$

\begin{theorem}\label{theoremB}
Let $K:{\mathcal D}\rightarrow\R^+$, $0<q<p<+\infty$, and $1<p<+\infty$.  Let $\mu$ and $\sigma$ be 
 nonnegative Borel measures on $\R^n$. Suppose that  $(K, \sigma) \in   {\rm (DLBO)}$. 
 Then there exists a constant $C>0$ such that the trace inequality
\begin{equation}\label{00.1}
\int_{\R^n} \left \vert \, T_{K_{\mathcal D}}[fd\sigma] \, \right \vert^q 
\,d\mu \leq  C \, \|f\|_{L^p(d\sigma)}^q,  \qquad f\in
L^p(d\sigma),
 \end{equation}
 holds if and only if 
\begin{equation}\label{00.2}
{\mathcal W}_{K, \,\sigma}^{\mathcal D} [\mu]\in
L^{\frac{q(p-1)}{p-q}}(d\mu).
 \end{equation}
Moreover, the least constant $C$ in {\rm (\ref{00.1})} is equivalent to 
$|| {\mathcal W}_{K, \,\sigma}^{\mathcal D} [\mu]||^{\frac q {p'}}_{L^{\frac{q(p-1)}{p-q}}(d\mu)}$ 
with constants of equivalence that depend only on $p$, $q$, and $A$.  
\end{theorem} 

\noindent {\bf Remark.} Note that the exponent $\frac{q(p-1)}{p-q}< 1$ in (\ref{00.2})   if $0<q<1$. 
In the case $q=1$, Theorem~\ref{theoremB} is equivalent by duality to Wolff's inequality {\rm (\ref{wolffinequality})}.
\bigskip

 We first give a characterization of the trace inequality in the diagonal case $q=p$ 
 in terms of discrete Carleson measures. It is deduced from the generalized Wolff 
 inequality and the dyadic Carleson measure theorem. For Riesz potentials, this characterization was 
  obtained earlier in  \cite{verbitsky2}.

\begin{lemma}\label{lemmaC}
Let $K:{\mathcal D}\rightarrow\R^+$. Let $1<p<+\infty$, 
 and let $\mu$ and $\sigma$ be nonnegative Borel
measures on $\R^n$. Suppose that  $(K, \sigma) \in   {\rm (DLBO)}$. 
Then the trace inequality  {\rm (\ref{00.1})} holds in the case $q=p$ 
if and only if  there exists a constant $B$ such that, for every dyadic cube $P \in {\mathcal D}$,  
 \begin{equation}\label{00.5}
 \sum_{Q \subset P} \,  K(Q) \, [\overline{K}(Q)]^{p'-1} \, [ \mu(Q)]^{p'} \, \sigma(Q) \le B \, \mu(P). 
\end{equation}
  
 \end{lemma}

 \begin{proof} By duality, {\rm(\ref{00.1})} holds with $q=p$ if and only if 
 \begin{equation}\label{00.6}
\int_{\R^n} \left \vert \, T_{K_{\mathcal D}}[g d\mu] \, \right \vert^{p'} (x)
\,d\sigma (x)\leq  C \, \|g\|_{L^{p'}(d \mu)}^{p'},  \qquad g\in
L^{p'}(d\mu).
\end{equation}
Without loss of generality we may assume that $g \ge 0$. By the generalized Wolff inequality, the left-hand side is 
equivalent to  
 \begin{equation}\label{00.7}
 \sum_{Q \in \D} \,  K(Q) \, [\overline{K}(Q)]^{p'-1} [\mu(Q)]^{p'} \sigma(Q) \, \left( \frac 1 {\mu(Q)} 
 \int_Q g \, d \mu\right)^{p'}. 
 \end{equation}
 Applying the dyadic Carleson measure theorem in $L^{p'} (d\mu)$ 
 (see, e.g., \cite{nazarovtreilvolberg} or \cite{verbitsky1}, Sec. 3), 
 we conclude that   (\ref{00.6}) 
 holds if and only if (\ref{00.5}) is valid. 
 \end{proof}

 \begin{lemma}\label{lemmaD}
 Let $1<p<+\infty$. Let $\mu$ and $\sigma$ be  nonnegative Borel measures on $\R^n$. 
  Suppose that  $(K, \sigma) \in   {\rm (DLBO)}$ with constant $A$ in (\ref{dlbo}). 
  Define $$d\mu_1=\frac{d\mu}{\left ( {\mathcal W}_{K, \,\sigma}^{\mathcal D} [\mu]\right)^{p-1}}.$$ Then
 \begin{equation}\label{lemmad0}
  \int_{\R^n} \left\vert T_{K_{\mathcal D}} [fd\sigma]\right \vert^p \,d\mu_1 \leq C \, \|f\|_{L^p(d\sigma)}^p, 
  \qquad f \in L^p(d\sigma), 
 \end{equation}
 where $C$ depends only on $p$ and $A$. 
 \end{lemma}

 \begin{proof}[Proof of Lemma~\ref{lemmaD}] By Lemma~\ref{lemmaC}, the preceding inequality 
 holds if and only if 
 (\ref{00.5}) is valid with $\mu_1$ in place of $\mu$. By H\"older's inequality,
 $$[\mu_1(Q)]^{p'} \le 
 [\mu(Q)]^{p'-1} \, \int_Q \frac{d \mu} {\left ( {\mathcal W}_{K, \,\sigma}^{\mathcal D} [\mu]\right)^{p}}.
 $$
 Then  by the preceding inequality, for any dyadic cube $P\in {\mathcal D}$,
 \begin{align}
 & \sum_{Q \subset P} \,  K(Q) \, [\overline{K}(Q)]^{p'-1} \, [ \mu_1(Q)]^{p'} \, \sigma(Q) \notag \\ 
 & \le  \sum_{Q \subset P} \,  K(Q) \, [\overline{K}(Q)]^{p'-1} \, [ \mu(Q)]^{p'-1} \, \sigma(Q) \notag \\ 
 & \times \int_Q \frac{d \mu_1 (y)} { \sum_{Q' \in {\mathcal D}} \, K(Q') \, 
 [\overline{K}(Q')]^{p'-1} \, [ \mu(Q')]^{p'-1} \, \sigma(Q') \, \chi_{Q'}(y)} \notag \\ & \le \mu_1(P).
 \notag
 \end{align}
 This proves (\ref{00.5}) with $\mu_1$ in place of $\mu$, and $B=1$. 
 \end{proof}
 
  We now prove Theorem~\ref{theoremB}.

 \begin{proof}[Proof of Theorem~\ref{theoremB}] We first prove that (\ref{00.2}) is sufficient in order 
 that (\ref{00.1}) hold. 
 H\"older's inequality with exponent $\frac{p}{q}>1$, together with Lemma~\ref{lemmaC}, yield 
  \begin{equation*}\begin{split}&
\left( \int_{\R^n} \left \vert  T_{K_{\mathcal D}} [fd\sigma]\right \vert^q \, d\mu \right)^\frac1{q}\\& \leq 
\left( \int_{\R^n} \left \vert  T_{K_{\mathcal D}}  [fd\sigma]\right \vert^p \frac{d\mu}{ {\mathcal W}_{K, \,\sigma}^d[\mu]^{p-1}}\right)^\frac1{p} 
\left(\int_{\R^n}   {\mathcal W}_{K, \,\sigma}^d[\mu]^\frac{q(p-1)}{p-q}d\mu \right)^\frac{p-q}{qp} \\& \leq C \,  \left( \int_{\R^n}  
 \left \vert T_{K_{\mathcal D}}[fd\sigma]\right \vert^p  d\mu_1\right)^\frac1{p} \leq C \, \|f\|_{L^p(d\sigma)}.
\end{split}\end{equation*}
 This proves (\ref{00.1}).

 We now prove the necessity of (\ref{00.2}). Suppose that (\ref{00.1}) holds. 
 For a sequence  of nonnegative reals $( \lambda_Q)_{Q \in {\mathcal D}}$ we set 
 $$
 f(x) =\sup_{Q\in{\mathcal D}} \, \lambda_Q  \, \chi_Q(x), \qquad x \in \R^n. 
 $$
 Obviously,  
$$ T_{K_{\mathcal D}} [fd\sigma] (x) \ge \sum_{Q \in {\mathcal D}} \, K(Q) \, \lambda_Q \, \sigma(Q) \, \chi_{Q} (x).
$$
Similarly, letting 
$$
 g(x)=\sum_{Q\in{\mathcal D}} \, \lambda_Q  \chi_Q(x),\qquad x \in \R^n, 
 $$
 we obtain
  \begin{align}
  T_{K_{\mathcal D}} [gd\sigma] (x) & \ge \sum_{Q \in {\mathcal D}} \, K(Q) \, \chi_Q(x) 
 \int_Q \left [ \sum_{Q' \supset Q} \, 
 \lambda_{Q'} \, \chi_{Q'}(y) \right ] \, d \sigma(y) \notag \\ & \ge \sum_{Q' \in {\mathcal D}} \lambda_{Q'} \, 
 \sigma(Q') \overline{K}(Q') \, \chi_{Q'}(x).
 \notag
 \end{align}
 Hence,  the following two estimates hold:
 \begin{equation}\begin{split}\label{estimate}&
 \left \Vert \sum_{Q \in {\mathcal D}}   \, K(Q) \, \lambda_Q \, \sigma(Q) \,
  \chi_{Q} \right \Vert_{L^q( d \mu)}   \le 
 C \, \left \Vert \sup_{Q\in{\mathcal D}} \, \lambda_Q  \chi_Q   \right \Vert_{L^p(d\sigma)},  \\ 
 &\left \Vert \sum_{Q \in {\mathcal D}}   \, 
 {\overline K} (Q) \, \lambda_Q \, \sigma(Q) \, \chi_{Q} \right \Vert_{L^q( d \mu)}  \le 
 C \, \left \Vert \sum_{Q\in{\mathcal D}} \, \lambda_Q  \chi_Q   \right \Vert_{L^p(d\sigma)}. 
 \end{split}\end{equation}
 
 Next, we rewrite the above two inequalities in terms of the discrete Littlewood--Paley 
 spaces ${\bf f}_{p\,\nu}^{0\,q}$, $0<p<+\infty$, $0<q\leq +\infty$, i.e., 
 spaces which are formed by the collection of all sequences $(s_Q)_{Q\in {\mathcal D}}$ such that 
 $$\Vert (s_Q)_Q\Vert_{{\bf f}_{p\,\nu}^{0\,q}}= 
 \left \Vert \left ( \sum_{Q\in {\mathcal D}_\mu} 
 (|s_Q|\nu(Q)^{-\frac12}\chi_Q)^q \right )^\frac1{q}\right \Vert_{L^p(d\nu)}
 <+\infty,$$
 where ${\mathcal D}_\nu$ is   the set of dyadic cubes $Q \in {\mathcal D}$ such that $\nu(Q)\neq 0$. 
 (See \cite{frazierjawerth}, \cite{verbitsky1}.) 
 
 Replacing $\lambda_Q$ by $\lambda_Q^s$, where  $s q>1$,  we obtain 
$$\left \Vert \sum_{Q \in {\mathcal D}}   \, K(Q) \, \lambda_Q^s \, \sigma(Q) \,
   \chi_Q\right \Vert_{L^q( d \mu)} =\left \Vert \left(K(Q)^\frac1{s} \lambda_Q \sigma(Q)^\frac1{s} \mu(Q)^\frac12\right)_{Q \in {\mathcal D}}  \right \Vert_{{\bf f}_{qs\,\mu}^{0\,s}}^s,$$ 
 $$\left \Vert \sum_{Q \in {\mathcal D}}   \, {\overline K}(Q) \, \lambda_Q^s \, \sigma(Q) \,
  \chi_{Q} \right \Vert_{L^q( d \mu)} =\left \Vert \left({\overline K}(Q)^\frac1{s} \lambda_Q \sigma(Q)^\frac1{s} \chi_Q \right)_{Q \in {\mathcal D}} \right \Vert_{{\bf f}_{qs\,\mu}^{0\,s}},$$
  and 
  $$\Vert \sup_{Q\in{\mathcal D}} \, \lambda_Q^s  \chi_Q    \Vert_{L^p(d\sigma)}=\Vert (\lambda_Q \sigma(Q)^\frac12)\Vert_{{\bf f}_{ps\,\sigma}^{0\,\infty}}^s,$$
  $$\Vert \sum_{Q\in{\mathcal D}} \, \lambda_Q^s  \chi_Q    \Vert_{L^p(d\sigma)}=\Vert (\lambda_Q \sigma(Q)^\frac12)\Vert_{{\bf f}_{ps\,\sigma}^{0\,s}}^s.$$  
  Consequently, the estimates in (\ref{estimate}) are equivalent to 
   \begin{align}
 \left \Vert \left(K(Q)^\frac1{s} \lambda_Q \sigma(Q)^{\frac1{s}-\frac12} \mu(Q)^\frac12 \right)_{Q \in {\mathcal D}}\right \Vert_{{\bf f}_{qs\,\mu}^{0\,s}}  & \le 
 C \, \left \Vert \left( \, \lambda_Q     \right)_{Q \in {\mathcal D}}\right \Vert_{{\bf f}_{ps\,\sigma}^{0\,\infty}}, \notag \\ 
 \left \Vert \left({\overline K}(Q)^\frac1{s} \lambda_Q \sigma(Q)^\frac1{s} \chi_Q \right)_{Q \in {\mathcal D}}\right \Vert_{{\bf f}_{qs\,\mu}^{0\,s}} & \le 
 C \, \left \Vert \left( \, \lambda_Q    \right)_{Q\in{\mathcal D}} \right \Vert_{{\bf f}_{ps\,\sigma}^{0\,s}}. \notag
 \end{align}

 It is shown in \cite{frazierjawerth}, Corollary 8.3, that if $1\leq p_0,q_0<+\infty$, $1\leq p_1,q_1\leq +\infty$, $0<\theta<1$, $1/p=(1-\theta)/p_0+\theta/p_1$ and $1/q=(1-\theta)/q_0+\theta/q_1$, then the complex interpolation  space $[{\bf f}_{p_0\,dx}^{0\,q_0},{\bf f}_{p_1\,dx}^{0\,q_1}]_\theta\simeq {\bf f}_{p \,dx}^{0\,q}$, where $dx$ is Lebesgue measure on $\R^n$. The main tool in their proof is the fact that the dyadic Hardy-Littlewood maximal operator (with respect to $dx$) is bounded in $L^p(dx)$, provided $p>1$. Since the  dyadic Hardy-Littlewood maximal operator with respect to an arbitrary nonnegative measure $d \mu$  is also bounded in $L^p(d \mu)$, it follows  that  $[{\bf f}_{p_0\,\sigma}^{0\,q_0},{\bf f}_{p_1\,\sigma}^{0\,q_1}]_\theta\simeq {\bf f}_{p \,\sigma}^{0\,q}$. (See 
 also \cite{cohnverbitsky}.) 
 
 Using  interpolation 
 for analytic families of operators acting on the pair of normed spaces ${\bf f}_{ps\,\sigma}^{0\,\infty}$ and ${\bf f}_{ps\,\sigma}^{0\,s}$ 
 (see, e.g.,  \cite{reedsimon},  Theorem IX.20) we deduce:
 $$\left \Vert \left(K(Q)^{\frac1{s}(1-\frac1{p})} {\overline K}(Q)^{\frac1{sp}} \lambda_Q \sigma(Q)^{\frac1{s}-\frac12}\mu(Q)^\frac12 \right)_{Q \in {\mathcal D}}\right \Vert_{{\bf f}_{qs\,\mu}^{0\,s} }   \le 
 C \, \left \Vert \left( \, \lambda_Q      \right)_{Q \in {\mathcal D}}\right \Vert_{{\bf f}_{ps\,\sigma}^{0\,ps}},$$
 and consequently 
 \begin{align} 
& \left \Vert \sum_{Q \in {\mathcal D}}   \, [K(Q)]^{1 - \frac 1 p} \,   [{\overline K} (Q)]^{\frac 1 p}  
 \, \lambda_Q \, \sigma(Q) \, \chi_{Q} \right \Vert_{L^q( d \mu)} \notag \\  & \le 
 C \, \left \Vert  \left ( \sum_{Q\in{\mathcal D}} \, \lambda_Q^p  \chi_Q \right )^{\frac 1 p}  
 \right \Vert_{L^p(d\sigma)} 
 = C \, \left [\sum_{Q\in{\mathcal D}} \, \lambda_Q^p \, \sigma(Q) \right ]^{\frac 1 p}.
 \notag
 \end{align}
 Now applying Theorem 3(c) in \cite{verbitsky1} we conclude:
 $$
 \sum_{Q\in{\mathcal D}} \,  [{\overline K} (Q)]^{p'-1}  \, K(Q) \, [\mu(Q)]^{p'-1} \, \sigma(Q) \, \chi_Q (x) 
= {\mathcal W}_{K, \,\sigma}^d[\mu] (x)    \in L^\frac{q(p-1)}{p-q}(d\mu).
 $$
 This completes the proof of Theorem~\ref{theoremB}. \end{proof}

 \section{Continuous trace inequalities}\label{section3}

  We recall that if $k:(0,+\infty)\rightarrow\R^+$ is a  nonincreasing lower
semicontinuous function, and   $\sigma$ is a nonnegative locally finite Borel
measure on $\R^n$, we have defined in the Introduction the function
 $$
\overline{k}(r)(x)=\frac1{\sigma(B(x,r))}\int_0^r
k(l) \, \sigma(B(x,l)) \, \frac{dl}{l},
$$ for $x\in\R^n$, $r>0$.
The continuous Wolff-type potential introduced in \cite{cascanteortegaverbitsky2} is given by: 
 \begin{equation}\label{wolffcontinuous}
 {\mathcal W}_{k,\,\sigma}[\mu](x)=\int_0^{+\infty} k(r)\sigma(B(x,r))
 \left(\int_{B(x,r)} \overline{k}(r)(y)d\mu(y)\right)^{p'-1}\frac{dr}{r}.
 \end{equation}
 
We will assume that $k$ and
 $\sigma$ satisfy   
 \begin{equation}\label{lbo}
  \sup_{y \in B(x,r)} \overline{k}(r)(y) \le A \, \inf_{y \in B(x,r)}
\overline{k}(r)(y),
 \end{equation}
 where $A$ does not depend on $x\in \R^n$ and $r>0$. In this case we will say that the pair $(k,\sigma)$
has the logarithmic bounded oscillation property, or  
simply write $(k,\sigma) \in \text{LBO}$.

  If $(k,\sigma)\in\text{LBO}$,   the Wolff-type potential can be defined in an equivalent form given by
 $$
 {\mathcal W}_{k,\,\sigma}[\mu](x)=\int_0^{+\infty} k(r)\sigma(B(x,r))
\overline{k}(r)(x)^{p'-1} \mu (B(x,r))^{p'-1}\frac{dr}{r}.
 $$
 
 If $\mu$ and $\sigma$ are nonnegative Borel measures on $\R^n$, and
$1<p<+\infty$, the energy of $\mu$ associated with $k$ and $\sigma$ is given by 
\begin{equation}\label{generalaenergy}  {\mathcal E}_{k,\,
\sigma}[\mu]=\int_{\R^n}\left( T_k[\mu](x)\right)^{p'}d\sigma(x). 
\end{equation} 
 
 We first observe that by setting $K(Q)=k(r_Q)$, we can associate to the radial kernel $k$ a dyadic kernel 
 $K_{\mathcal D}$, 
 and the corresponding integral operator $T_{K_{\mathcal D}}$. The following proposition gives   
 a relationship between the continuous energy and a supremum of the dyadic energy over shifts 
 $ {\mathcal D}+z$, $z\in\R^n$, where the shifted dyadic Wolff type potential is defined by
  \begin{align}\label{formula1.2.1}
{\mathcal{W}}_{K, \,\sigma}^{{\mathcal D}+z}[\mu](x) =\sum_{Q\in{\mathcal D}} & K(Q+z)\sigma(Q+z) \notag \\ & \times \left(   \int_{ Q+z}
\overline{K}(Q+z)(y)d\mu(y)  \right)^{p'-1}\chi_{Q+z}(x), \quad x\in\R^n.
\notag
\end{align}
\begin{proposition}\label{propositionA}
Let  $k:(0,+\infty)\rightarrow\R^+$
be a  nonincreasing lower semicontinuous function. Let 
$1<p<+\infty$, and  let $\sigma$ be  a nonnegative
locally finite Borel measure on $\R^n$. Suppose that $\sigma$  satisfies a 
doubling condition and that  $(k,\sigma)\in${\rm\text{LBO}}.  If $K(Q)=k(r_Q)$,  then for any
nonnegative Borel measure $\mu$ on $\R^n$, 
\begin{equation}\label{wolffgeneral}{\mathcal E}_{k, \, \sigma}[\mu]\simeq
\sup_{z\in\R^n}\int_{\R^n}{\mathcal W}_{K, \, \sigma}^{{\mathcal D}+z}[\mu](x) \, d\mu(x),\end{equation} with
constants of equivalence that may depend on $k$ and $\sigma$, but not on $\mu$. 

\end{proposition}
  
  \begin{proof}[Proof of Proposition \ref{propositionA}] In \cite{cascanteortegaverbitsky2}, Proposition 3.7, it is proved that there exist constants
   $c, \, C>0$ such that for any
  $x\in\R^n$,  
  \begin{equation}\label{01.1}\sum_{Q\in{\mathcal D}}k(cr_Q)\sigma(Q)\chi_Q(x)\left(
\int_Q\overline{ K }(Q)(y)d\mu(y)\right)^{p'-1}\leq C \, {\mathcal W}_{k,\,\sigma}[\mu](x). \end{equation} 
  The constant $c$ above appears because, if $Q$ is a cube and $x\in Q$, then 
  $Q\subset B(x,\frac{c}2r_Q)$. Since $\overline K$ satisfies a doubling condition, the constant $c$ is 
  not needed in $\overline{K}$.
  
   We also have (Lemma 3.9 in \cite{cascanteortegaverbitsky2}) that for any $c>0$,
  \begin{equation}\label{pointwise}\sum_{Q\in{\mathcal
D}}k(cr_Q)\sigma(Q)\overline{k}(r_Q)^{p'-1}\mu(Q)^{p'}\simeq 
\sum_{Q\in{\mathcal D}}k(r_Q)\sigma(Q)\overline{k}(r_Q)^{p'-1}\mu(Q)^{p'},\end{equation}
with constants that do not depend on $\mu$.
 Consequently, (\ref{01.1}) and (\ref{pointwise}) give that if we write $k_c(r)=k(cr)$, and 
 $K_c(Q)=k_c(r_Q)$, then 
 \begin{equation*}\begin{split} {\mathcal E}_{k, \, \sigma}[\mu] & \simeq \int_{\R^n}{\mathcal W}_{k,\,\sigma}[\mu](x)d\mu(x)  \\&\geq C\sup_{z\in \R^n} \sum_{Q\in{\mathcal D}}K_c({Q+z})\sigma(Q+z) \overline{ K_c }(Q+z)^{p'-1}\mu(Q+z)^{p'-1}\\ & \simeq  
 \sup_{z\in\R^n}\int_{\R^n}{\mathcal W}_{K,\,\sigma}^{{\mathcal D}+z}[\mu](x)d\mu(x).\end{split}\end{equation*}
 
 On the other hand, the estimate obtained in  \cite{cascanteortegaverbitsky2}, 
 page 870, together with (\ref{pointwise}), gives
 $${\mathcal E}_{k, \, \sigma}[\mu]=\int_{\R^n} T_k[\mu](x)^{p'}d\sigma(x)
  \leq C\sup_{z\in\R^n} \sum_{Q\in{\mathcal D}} K(Q+z)\sigma(Q+z) 
 \overline{ K }(Q+z)^{p'-1}\mu(Q+z)^{p'}. $$
 \end{proof}

 Our next lemma shows that  
  if the trace inequality holds for $T_k$, and $0<q<+\infty$, $1<p<+\infty$, then it is 
  also valid for any operator $T_{k_c}$ associated with the kernel $k_c(r)=k(cr)$, $c>0$. 
 
 \begin{lemma}\label{lemmatilde} Let  $k:(0,+\infty)\rightarrow\R^+$
be a  nonincreasing lower semicontinuous function, and let $k_c(x) = k(c x)$, where $c>0$. Let 
$1<p<+\infty$ and  $0<q< \infty$. Suppose that $\sigma$ and $\mu$ are  nonnegative Borel 
measures on $\R^n$, and  $\sigma$ is a  doubling measure. 
  Then the following statements are  equivalent:
  
  (i) The inequality 
  \begin{equation} \label{trace1} 
  \int_{\R^n} |T_k [fd\sigma]|^q \,d\mu\leq C_1 \, \|f\|_{L^p(d\sigma)}^q
 \end{equation}
 holds for all $f \in L^p(d\sigma)$. 
 
 (ii) The inequality  
  \begin{equation} \label{trace2} 
  \int_{\R^n} |T_{k_c} [fd\sigma]|^q \, d\mu \leq C_2 \, \|f\|_{L^p(d\sigma)}^q
 \end{equation}
 holds for all $f \in L^p(d\sigma)$.

 Moreover, the least constants  $C_i$, $i=1,2$, in the above inequalities are equivalent, and the constants of equivalence depend only on $k$, $n$, $q$, 
 $p$, $c$, and the doubling constant of $\sigma$. \end{lemma}
 
\begin{proof}[Proof of Lemma \ref{lemmatilde}] Since $k$ is nonincreasing it suffices to prove  that (\ref{trace1}) implies (\ref{trace2}) for $c$ small enough.
  Without loss of generality we may assume that 
 $f \ge 0$ and $c= \frac 1 2$.  
 Denote by $M^\sigma$ the centered maximal function  with respect to $\sigma$ on $\R^n$ defined by 
$$M^\sigma f(x) = \sup_{r>0} \, \frac 1 {\sigma(B(x,r))} \int_{B(x,r)} |f| \, d \sigma.$$
  We then have 
 \begin{align}
 T_k [(M^\sigma f) d \sigma] (x) & \ge \sum_{l \in \mathbb Z} \, \int_{2^{l-2} < |x-y| \le 2^{l-1}} k(x-y) M^\sigma f(y) \, 
 d \sigma (y) \notag \\ & \ge \sum_{l \in \mathbb Z} \, k(2^{l-1}) \int_{2^{l-2} < |x-y| \le 2^{l-1}}  M^\sigma f(y) \, 
 d \sigma (y).\notag
 \end{align}
 Notice that, for $y \in B(x, 2^{l-1})$ and $r =2^{l+2}$, 
 $$ B(x, 2^{l+1}) \subset B(y, r) \subset B(x, 2^{l+3}).$$
 Hence, for   $ d \nu = f \, d \sigma$, we have 
 $$M^\sigma f (y) \ge \frac 1 {\sigma(B(y,r))} \int_{B(x,2^{l+1})} f \, d \sigma \ge C \, 
 \frac 1 {\sigma(B(x, 2^{l+3}))} \, \nu(B(x, 2^{l+1})).$$
 Since $\sigma$ is a doubling measure, it follows (see \cite{cascanteortegaverbitsky2}, p. 874) that 
 $$\sigma ( \{y: \, 2^{l-2} < |x-y| \le 2^{l-1} \}) \simeq \sigma (B(x, 2^{l+3})).$$
  Thus, 
 $$T_k [(M^\sigma f) d \sigma] (x) \ge C \, \sum_{l \in \mathbb Z} \,k(2^{l-1}) \, \nu(B(x, 2^{l+1})).$$
 It is well-known (see, e.g., \cite{fefferman}) that $M^\sigma$ is a bounded operator on $L^p(d \sigma)$, $p>1$. Thus, 
 it follows that  (\ref{trace1}) implies 
  \begin{equation} \label{trace4}
   \int_{\R^n} \left \vert \sum_{l \in \mathbb Z} \,  k(2^{l-1}) \int_{B(x, \, 2^{l+1})} 
   f \, d \sigma \right \vert^q \, d\mu \leq C \, \|f\|_{L^p(d\sigma)}^q.
 \end{equation}
  It remains to check that in its turn (\ref{trace4}) implies (\ref{trace2}). Indeed, this follows 
  from the estimate 
  $$T_{k_\frac12} [f d \sigma] (x) = \sum_{l \in \mathbb Z} \, \int_{2^{l} < |x-y| \le 2^{l+1}} k_\frac12(x-y)
  \, f(y) \, d \sigma(y)  \le \sum_{l \in \mathbb Z} \, k_\frac12(2^{l}) \, \nu(B(x, 2^{l+1})).$$
 The proof of Lemma~\ref{lemmatilde} is complete. \end{proof}

 \begin{theorem}\label{theoremC}
Let $K:{\mathcal D}\rightarrow\R^+$, $0<q<p<+\infty$, and $1<p<+\infty$.  Let $\mu$ and $\sigma$ be 
nonnegative Borel measures on $\R^n$. Suppose that  $(k, \sigma) \in   {\rm (LBO)}$, and $\sigma$ 
satisfies a doubling condition. 
 Then the trace inequality
\begin{equation}\label{00.3}
\int_{\R^n} \left \vert \, T_{k}[fd\sigma] \, \right \vert^q 
\,d\mu \leq  C \, \|f\|_{L^p(d\sigma)}^q,  \qquad f\in
L^p(d\sigma), 
 \end{equation}
 holds if and only if   $ \int_{\R^n}({\mathcal W}_{k, \,\sigma}[\mu])^{\frac{q(p-1)}{p-q}}d\mu<+\infty$.
 \end{theorem} 
  
 \begin{proof}[Proof of Theorem \ref{theoremC}] The proof of the sufficiency is immediate from the following lemma, which is a continuous 
 version of Lemma \ref{lemmaC}.
 \begin{lemma}\label{lemmaF} Under the same hypotheses as in  Theorem \ref{theoremC}, if 
 $$d\mu_1=\frac{d\mu}{( {\mathcal W}_{k, \,\sigma}[\mu])^{p-1}},$$ then
 \begin{equation}\label{lemmapp1}
   \int_{\R^n} | T_k[fd\sigma]|^p (x)\,d\mu_1(x)\leq C \, \|f\|_{L^p(d\sigma)}^p.
 \end{equation}
 \end{lemma}

\begin{proof}[Proof of Lemma \ref{lemmaF}] Duality and Lemma \ref{lemmatilde} show that  (\ref{lemmapp1}) holds if and only if 
$${\mathcal E}_{k_c, \,\sigma}[gd\mu_1]\leq C \|g\|_{L^{p'}(d\mu_1)}^{p'},$$ 
where $k_c(r)=k(cr)$, and $c>0$ is the constant which appears in (\ref{01.1}).
  Proposition  \ref{propositionA}   applied to $g\mu_1$ and $k_c$ gives 
$${\mathcal E}_{k_c, \, \sigma}[g\mu_1]\simeq
\sup_{z\in\R^n}\int_{\R^n}{\mathcal W}_{K_c, \, \sigma}^{{\mathcal D}+z}[g\mu_1](x) \, g(x)d\mu_1(x)\simeq \sup_{z\in\R^n} {\mathcal E}_{K_c,\sigma}^{{\mathcal D}+z}[g\mu_1].$$

Consequently, it remains to show that 
\begin{equation*}\label{lemmapp2}
   \sup_{z\in\R^n}\int_{\R^n}{T_{K_c}}_{{\mathcal D}+z}[fd\sigma]^p (x)\,d\mu_1(x)\leq C\|f\|_{L^p(d\sigma)}^p.
 \end{equation*}
But (\ref{pointwise})  yields that, for each $z\in\R^n$, $d\mu_1(x)\leq C \frac{d\mu(x)}{{\mathcal W}_{{K_c}, \, \sigma}^{{\mathcal D}+z}[\mu](x)^{p-1}}$, with constant that does not depend on $z$.
  Applying Lemma \ref{lemmaD}, we complete the proof of Lemma~\ref{lemmaF}.
  \end{proof}

 Now the proof of the sufficiency part in Theorem \ref{theoremC} follows the same argument as in 
 the case of the discrete trace inequality. 
 
 To establish the necessity part, we start with a lemma. 
 \begin{lemma}\label{lemmadiscrete}
 Assume that (\ref{00.3}) holds. Then there exist $c>0$ which depends only on the dimension, 
 and $C>0$ such that,  for any sequence $(\lambda_Q)$ of nonnegative reals, and any $z\in\R^n$, 
 \begin{align}\label{interpol}
& \left\Vert \sum_{Q\in{\mathcal D}}[K_c(Q+z)]^{1-\frac 1{p}} [\overline{K_c}(Q+z)]^{\frac 1{p}}  \lambda_Q \,  \sigma(Q+z) \, \chi_{Q+z} \right\|_{L^q(d\mu)} \notag \\ & \leq C \left[\sum_{Q\in {\mathcal D}} \lambda_Q^p \sigma(Q+z)\right]^{\frac1{p}}.
  \end{align}
   \end{lemma}
 \begin{proof}[Proof of Lemma \ref{lemmadiscrete}] If $(\lambda_Q)_Q$ is a sequence of nonnegative reals, we define 
 $f(x)= \sup_{Q\in{\mathcal D}} \lambda_Q\chi_Q(x)$. 
 For any $i\in\Z$, we denote by $Q_i$ the unique cube in ${\mathcal D}$ suich 
 that $x\in Q_i$ and $|Q_i|=2^{in}$. We then have
  \begin{equation*}\begin{split}& T_kf(x)=\int_{\R^n} k(|x-y|)\sup_{Q\in {\mathcal D}} \lambda_Q\chi_Q(y)d\sigma(y)\geq  \\&
 \sum_{i\in\Z} \lambda_{Q_{i+1}} \int_{Q_{i+1}\setminus Q_i} k(|x-y|) \chi_{Q_{i+1}}(y)d\sigma(y)= \sum_{i\in\Z} \lambda_{Q_{i+1}} \int_{Q_{i+1}\setminus Q_i} k(|x-y|) d\sigma(y).
 \end{split}\end{equation*}
 Next, there exists $c>0$ depending only on $n$, such that, for any $y\in Q_{i+1}$, $|x-y|\leq c2^{i+1}=cr_{Q_{i+1}}$. Since $\sigma$ satisfies a doubling condition, we also have that $\sigma(Q_{i+1}\setminus Q_i)\simeq \sigma(Q_{i+1}).$ Altogether, we deduce that the above sum is bounded from below by
 $$C\sum_{i\in\Z} \lambda_{Q_{i+1}}k_c(r_{Q_{i+1}}) \sigma(Q_{i+1})\chi_{Q_{i+1}}(x)=\sum_{Q\in{\mathcal D}} \lambda_Q K_c(Q)\sigma(Q)\chi_Q(x) .$$
 On the other hand, let $g(x)=\sum_{Q\in{\mathcal D}}\lambda_Q\chi_Q(x)$.  We then have
  $$ T_kg(x)=\int_{\R^n} k(|x-y|)\sum_{Q\in{\mathcal D}} \lambda_Q \chi_Q(y)d\sigma(y)\geq \sum_{x\in Q}\lambda_Q \int_Q k(|x-y|) d\sigma(y).$$
 If $x\in Q$, we denote by $2^{-l}Q$, $l\geq0$, the unique cube in ${\mathcal D}$ satisfying $x\in 2^{-l}Q$ and $r_{2^{-l}Q}=2^{-l}r_Q$. Then
  \begin{equation*}\begin{split}& \int_Qk(|x-y|)d\sigma(y) = \sum_{l\geq 0} \int_{2^{-l}Q\setminus 2^{-l-1}Q} k(|x-y|)d\sigma(y)\geq \\&
  \sum_{l\geq 0} k_c(2^{-l}r_Q) \sigma(2^{-l}Q\setminus 2^{-l-1}Q)\simeq 
 \sum_{l\geq 0} k_c(2^{-l}r_Q) \sigma(2^{-l}Q)\simeq \overline{K_c}(Q)\sigma(Q).
 \end{split}\end{equation*}
  Consequently, $$T_kg(x)\geq C \sum_{Q\in {\mathcal D}} \lambda_Q \sigma(Q) \overline{K_c}(Q)\chi_Q(x).$$
  We then have:
  \begin{equation*}\begin{split}& \left \Vert\sum_{Q\in{\mathcal D}} \lambda_Q K_c(Q) \, \sigma(Q)\chi_Q 
  \right \Vert_{L^q(d\mu)} 
  \leq C \, \left \Vert\sup_{Q\in{\mathcal D}} \lambda_Q\chi_Q\right\Vert_{L^p(d\sigma)}, \\ 
  &\left \Vert \sum_{Q\in {\mathcal D}}\lambda_Q\overline{K_c} (Q) \, 
  \sigma(Q)\chi_Q\right \Vert_{L^q(d\mu)} \leq C \, \left \Vert \sum_{Q\in{\mathcal D}} \lambda_Q\chi_Q\right \Vert_{L^p(d\sigma)}.
   \end{split}\end{equation*}

 We can now employ the interpolation argument used in the proof of 
 the necessity for the discrete trace inequality to obtain
  $$\left \Vert\sum_{Q\in{\mathcal D}}[K_c(Q)]^{1-\frac1{p}} [\overline{K_c}(Q)]^\frac1{p} 
  \lambda_Q \sigma(Q)\chi_Q\right \Vert_{L^q(d\mu)}\leq C\left[\sum_{Q\in {\mathcal D}} \lambda_Q^p 
  \sigma(Q)\right]^\frac1{p}.$$
  Applying the same argument to a shifted 
 dyadic lattice $\mathcal D + z$, $z \in \R^n$, we obtain a similar inequality with a constant $C$ 
 which does not depend on $z$. \end{proof}
 
 In order to follow with the proof of the necessity, we  recall that by Theorem 3 (c) in \cite{verbitsky1}, 
  estimate (\ref{interpol}) is equivalent to the condition 
  ${\mathcal W}_{K_c, \,\sigma}^{{\mathcal D}+z}[\mu]    \in L^\frac{q(p-1)}{p-q}(d\mu), 
 $ 
 uniformly in $z$, where 
$$
{\mathcal W}_{K_c, \,\sigma}^{{\mathcal D}+z}[\mu] (x) = \sum_{Q\in{\mathcal D}} \,  [{\overline {K_c}}(Q+z)]^{p'-1}  \, K_c(Q+z) \, [\mu(Q+z)]^{p'-1} \, \sigma(Q+z) \, \chi_{Q+z} (x).
 $$
 
 Our next goal is to show that the proof of the above result can be modified to obtain that 
 if $B_j = B(0, 2^{j+j_0})$ and $x \in B_j$ (here $j_0$ is fixed so that $2^{j_0} > 2 \sqrt{n} +1$), then
 \begin{equation}\label{averages}
\sum_{Q \in {\mathcal D}} \left [ \frac 1 {|B_j|}\int_{B_j} 
\left ( [{\overline K_c} (Q+z)]^{p'-1}  \, K_c(Q+z) \, [\mu(Q+z)]^{p'-1} \, \sigma(Q+z) \, 
\chi_{Q+z} (x) \, \right)^{\frac {p-1}{ps-1}} dz \right]^{\frac {ps-1} {p-1}} 
\end{equation}
belongs to $L^{\frac{q(p-1)}{p-q}} (d \mu)$, where $s>1$ satisfies $sq>1$.

Indeed, if $s>1$ is chosen so that $sq>1$, (\ref{interpol}) can be rewritten as
$$\sup_{z\in \R^n}\left \Vert \sum_{Q\in{\mathcal D}}[K_c(Q+z)]^{(1-\frac1{p})\frac1{s}} [\overline{K_c}(Q+z)]^{\frac1{ps}} \lambda_Q \sigma(Q+z)^\frac1{sp'}\chi_{Q+z}\right\Vert_{L^{qs}(l^{ps})(d\mu)}\leq C\| (\lambda_Q) \|_{l^{ps}}.$$
 We  then have, arguing  as  in \cite{verbitsky1}, p. 545, that the above inequality is equivalent to 
 \begin{equation*}\begin{split}
 & \sup_{z\in \R^n}\left(\sum_{Q\in{\mathcal D}}  
 \left([K_c(Q+z)]^{(1-\frac1{p})\frac1{s}}   
 [\overline{K_c}(Q+z)]^{\frac1{ps}}  \sigma(Q+z)^\frac1{sp'}\right)^\frac{ps}{ps-1} 
 \left \vert \int_{Q+z} g_Qd\mu\right \vert^\frac{ps}{ps-1}\right)^\frac{ps-1}{ps}\\& \leq 
 C \, ||g||_{L^\frac{qs}{qs-1}(l^{s'})(d\mu)},
 \end{split}
 \end{equation*}
 for all $g=(g_Q)\in L^\frac{qs}{qs-1}(l^{s'})(d \mu) $. Here $L^s(l^r)(d \mu)$ denotes the 
 mixed-norm Lebesgue space of vector-valued functions $g=(g_Q)$ equipped with the quasi-norm 
 $$||g||_{L^s(l^r)(d \mu)} = \left [ \int_{\R^n} ( \sum_Q |g_Q|^r)^{\frac s r} \, d \mu\right]^{\frac 1 s}.
 $$

 Setting $g_Q=\psi_Q^{\frac{ps-1}{ps}}$, $\psi_Q\geq 0$, and using H\"older's inequality together with the dyadic vector-valued Fefferman-Stein maximal 
 theorem, it is easy to see that the last estimate is equivalent to
 \begin{equation*}\begin{split} \sup_{z\in \R^n} \sum_{Q\in{\mathcal D}} & 
 \left([K_c(Q+z)]^{(1-\frac1{p})\frac1{s}}[\overline{K_c}(Q+z)]^{\frac1{ps}}  \sigma(Q+z)^\frac1{sp'}\mu(Q+z)\right)^\frac{ps}{ps-1}  \\ & \times \frac1{\mu(Q+z)} \int_{Q+z} \psi_Qd\mu\leq 
 C \, ||\psi||_{L^{\frac{ps-1}{qs-1}\frac{q}{p}} (l^{\frac{ps-1}{s-1}\frac1{p}})(d \mu) },
 \end{split}\end{equation*}
 for all $\psi=(\psi_Q)_Q\in L^{\frac{ps-1}{qs-1}\frac{q}{p}}_\mu(l^{\frac{ps-1}{s-1}\frac1{p}})$.
(See the proof of Theorem 3(c) in \cite{verbitsky1} for details.) Averaging over $z \in B_j$  
gives 
\begin{equation*}\begin{split}\sum_{Q\in{\mathcal D}} & 
\frac1{|B_j|}\int_{B_j}\left([K_c(Q+z)]^{(1-\frac1{p})\frac1{s}}[\overline{K_c}(Q+z)]^{\frac1{ps}}  \sigma(Q+z)^\frac1{sp'}\mu(Q+z)\right)^\frac{ps}{ps-1} \times \\ & \frac1{\mu(Q+z)} \int_{Q+z} \psi_{Q} 
\, d\mu \, dz\leq 
 C \, ||\psi||_{L^{\frac{ps-1}{qs-1}\frac{q}{p}} (l^{\frac{ps-1}{s-1}\frac1{p}})(d \mu) }.
 \end{split}\end{equation*}
The preceding estimate by duality yields (\ref{averages}).

Let $c>0$ be the constant  in Lemma \ref{lemmadiscrete}. 
As was shown in \cite{cascanteortegaverbitsky2}, p. 877,  
 for the truncated Wolff type potential associated with $k_{2c}$ defined by 
 $${\mathcal W}_{k_{2c},\,\sigma}^R[\mu](x)= \int_0^R k_{2c}(r)\sigma(B(x,r)){\overline{k_{2c}}}(r)^{p'-1}(x)\mu(B(x,r))^{p'-1}\frac{dr}{r},$$
it follows that, for every $l \le j$ and $x \in B_j=B(0,2^j)$,
$${\mathcal W}_{k_{2c},\,\sigma}^{2^j}[\mu](x)\leq C \sum_{l\leq j} \sigma(B(x,2^l))k_{2c}(2^{l-1}){\overline{k_{2c}}}(2^l)(x)^{p'-1} \mu(B(x,2^l))^{p'-1}.$$
Next, for any $l\leq j$,
\begin{align}
&  [{\overline k_{2c}} (2^{l})]^{p'-1}  \, k_{2c}(2^{l-1}) \, [\mu(B(x, 2^{l}))]^{p'-1} \, 
\sigma(B(x, 2^{l})) \preceq\notag \\ & \left( \frac 1 {|B_j|}\int_{B_j} \, 
\left (\,  [{\overline K_c} (Q_l+z )]^{p'-1}  \, K_c(Q_l+z ) \, [\mu(Q_l+z )]^{p'-1} \, \sigma(Q_l+z ) \, \chi_{Q_l+z} (x) \right)^\frac{p-1}{ps-1}\, dz \right)^\frac{ps-1}{p-1},
\notag
\end{align}
where $Q_l+z$ is the unique cube of side length $r(Q+z) = 2^{l+1}$  that contains $x$. This follows as in the proof of Theorem 3.13 in \cite{cascanteortegaverbitsky2} by noticing that   
 all sums over the generation of cubes $Q+z$ 
such that $r_{Q+z} = 2^{l+1}$ on p. 877 in  \cite{cascanteortegaverbitsky2}  actually contain only one term. Hence,
$$
\int_{B_j}  \left (  \sum_{l \le j} [{\overline k_{2c}} (2^{l})]^{p'-1}  \, k_{2c}(2^{l-1}) \, [\mu(B(x, 2^{l}))]^{p'-1} \, 
\sigma(B(x, 2^{l}))  
\right)^{\frac{q(p-1)}{p-q}} d \mu \le C < +\infty.
$$
Letting $j \to +\infty$, we obtain  $ \int_{\R^n}({\mathcal W}_{k_{2c}, \,\sigma}[\mu])^{\frac{q(p-1)}{p-q}}d\mu<+\infty$. 

Now using the same argument with the kernel $k_\frac1{2c}$ in place of $k$ and applying   Lemma \ref{lemmatilde},  we complete the proof of Theorem~\ref{theoremC}.\end{proof}

\end{document}